%


\documentstyle{amsppt}
\magnification=\magstep1
\NoBlackBoxes
\hsize 6.0truein
\vsize 8.2truein
\hcorrection{.2truein}
\vcorrection{.2truein}
\loadbold

\define\pis{$\Pi^1_2$-singleton}

\define\jiny{(i_1,\cdots,i_{n+1})}
\document
\baselineskip=15pt

\font\bigtenrm=cmr12 scaled\magstep2
\centerline
{\bigtenrm {Provable $\Pi^1_2$-Singletons}}

\vskip20pt

\font\bigtenrm=cmr10 scaled\magstep2
\centerline
{\bigtenrm{Sy D. Friedman}\footnote"*"{Research supported by NSF contract \#
9205530. }} 

\centerline
{\bigtenrm {M.I.T.}}

\vskip20pt

\comment
lasteqno 1@0
\endcomment

In Friedman [90] we constructed a $\Pi^1_2$-singleton $R$  such that
$0<_LR<_L0^{\#}.$ An open question is whether such a $\Pi^1_2$-singleton
can be $ZFC$-provable, in the sense that $ZFC\vdash\phi$ has at most one
solution, where
$\phi$ is a $\Pi^1_2$ formula characterizing $R.$  In this note
we observe that the construction from Friedman [90] can be used to obtain a
$T$-provable \pis \, $R,0<_LR<_L0^{\#},$ where $T$ is a theory consistent with
$V=L$ contained as a subtheory of $ZFC+0^{\#}$ exists. $ T$ has consistency
strength  approximately that of $ZFC+$ there exists an ineffable
cardinal.

First we recall a definition from Friedman [90]. For $i_1<\cdots<i_{n+1},n\ge
1$  define  $I(i_1,\cdots,i_{n+1})=\{i<i_1|i$ is $L$-inaccessible and $i,i_1$
satisfy the same $\Sigma_1$ properties in $L_{i_{n+1}}$ with parameters from
$i\cup\{i_2,\cdots, i_n\}\}.$  An {\it acceptable quess} is such a sequence
$(i_1,\cdots,i_{n+1})$  
where $i_1$ is $L$-inaccessible and $1\le k<\ell\le n\longrightarrow i_k\in
I(i_{\ell},\cdots,i_{n+1}).$ 

Now we say that an acceptable  guess $(i_1,\cdots,i_{n+1})$  
is {\it good} if in addition $I\jiny$ is stationary in $i_1.$ 
We refer to $n$  as the {\it length} of the guess $\jiny$.

$T$ is the theory $ZFC+$ There are arbitrarily long good guesses. $T$ is a
subtheory of $ZFC+0^{\#}$ exists since any increasing sequence of Silver
indiscernibles $\jiny,$ where $n\ge 1$  and $i_1$ is regular, is a good guess.
(In fact I$\jiny$ is $CUB$ in $i_1$ in this case.) Also note that $T$ follows
from the existence, for each $n,$  of a cardinal $K$  such that any 
function on
$n$-tuples from $K$  has a homogeneous set $X$  containing an $\alpha$ such
that $X\cap\alpha$ is stationary in $\alpha,$  together with $n-2$  larger
ordinals. (This is close in strength to an ineffable cardinal.) And if $T$  is
true then it is true in $L.$  

Now recall that in Friedman [90] a \pis\, $R$ is constructed so as to ``kill''
acceptable guesses $\jiny$ such that $i_{n+1}<(i^+)^L$ and $p\jiny_0$
contradicts $R.$  Here, $p\jiny$ is a $\Sigma_1(L)$-procedure which assigns a
forcing condition to the guess $\jiny$ and $p\jiny_0$ is the ``real part'' of
$p\jiny,$ consisting of a function from $(2^{<\omega })^{<\omega }$ into
perfect trees. $R$ is in fact  a set of finite sequences of finite sequences
of $0$'s and $1$'s and is determined by the $p\jiny_0$ where $p\jiny$ belongs
to the generic class. A simple requirement that we may impose on the procedure
$p\jiny$ is that $p\jiny$ must decide which of the first $n$ elements of
$(2^{<\omega })^{<\omega }$ belongs to $R,$  for some fixed (constructible)
$\omega$-listing of $(2^{<\omega })^{<\omega }.$ 

An acceptable guess $\jiny$ is killed by adding a $CUB$ subset to $i_1$ 
disjoint from $I\jiny$. The $\Pi^1_2$ formula characterizing $R$  implies that
$R$ kills all acceptable guesses $\jiny$ such that $i_{n+1}<(i_1^+)^L$ 
and $p\jiny$ forces a false membership fact about $R.$  Now suppose $T$  holds
and that $R\ne S$ were both solutions to our $\Pi^1_2$ formula. Choose $n$  so
that $R$  and $S$  differ on the membership of one of the first $n$  elements
of $(2^{<\omega })^{<\omega }$ and let $\jiny$ be a good guess. By a Skolem
hull argument we may assume that $i_{n+1}<(i_1^+)^L.$ Then either $R$ or $S$
must kill $\jiny$ since $p\jiny$ decides membership of the first $n$  elements
of $(2^{<\omega })^{<\omega }.$ But goodness means that $I\jiny$ 
is stationary, a contradiction.

So $T$  proves that our $\Pi^1_2$ formula characterizing $R$  has at most one
solution.

\vskip20pt

\centerline{\bf Reference}

\vskip5pt

Friedman [90] \qquad The $\Pi^1_2$-Singleton Conjecture, Journal of the AMS.

\enddocument